\documentclass[a4paper, 11pt]{article}

\usepackage[latin1]{inputenc}
\usepackage[english]{babel}
\usepackage{graphicx}
\usepackage{amssymb}
\usepackage{amsmath}
\usepackage{amsfonts}
\usepackage{amsthm}
\usepackage{url}

\newtheorem{algorithm}{Algorithm}

\newtheorem{remark}{Remark}
\numberwithin{equation}{section}

\usepackage{algorithm}
\usepackage{algorithmic}
\def \matlab    {MATLAB$^{\textrm{\tiny \textregistered}}$}

\title{A Multilevel Approach For Nonnegative Matrix Factorization}

\author{\normalsize Nicolas Gillis${}^1$ and Fran\c{c}ois Glineur${}^2$}

\date{}

\begin{document}

\maketitle

\begin{abstract}
Nonnegative Matrix Factorization (NMF) is the problem of approximating a nonnegative matrix with the product of two low-rank nonnegative matrices and has been shown to be particularly useful in many applications, \rm{e.g.}, in text mining, image processing, computational biology, etc. 
In this paper, we explain how algorithms for NMF can be embedded into the framework of multilevel methods in order to accelerate their convergence. 
This technique can be applied in situations where data admit a good approximate representation in a lower dimensional space through linear transformations preserving nonnegativity. 
A simple multilevel strategy is described and is experimentally shown to speed up significantly three popular NMF algorithms (alternating nonnegative least squares, multiplicative updates and hierarchical alternating least squares) on several standard image datasets.
\bigskip

\end{abstract}

\footnotetext[1] {University of Waterloo, Department of Combinatorics and Optimization, Waterloo, Ontario N2L 3G1, Canada. E-mail: ngillis@uwaterloo.ca. This work was carried out when the author was a Research fellow of the Fonds de la Recherche Scientifique (F.R.S.-FNRS) at Universit\'e catholique de Louvain.} 

\footnotetext[2] {Universit\'e catholique de Louvain, CORE and ICTEAM Institute, B-1348 Louvain-la-Neuve, Belgium. 
E-mail: francois.\mbox{glineur}@uclouvain.be.  This text presents research results of the Belgian Program on Interuniversity Poles of Attraction initiated by the Belgian State, Prime Minister's Office, Science Policy
Programming. The scientific responsibility is assumed by the authors.}

\section{Introduction} \label{intro}

\begin{figure*}[ht!]
\begin{center}
\includegraphics[width=\textwidth]{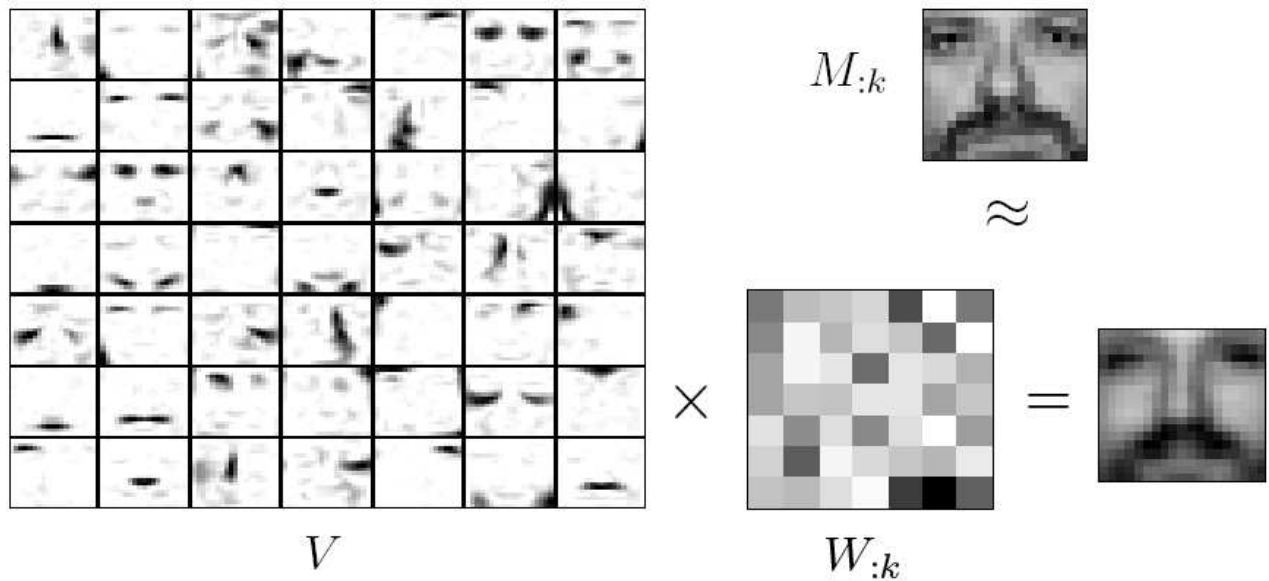}
\caption{Illustration of NMF on a face database. Basis elements (matrix $V$) obtained with NMF on the CBCL Face Database $\#1$, MIT Center For Biological and Computation Learning, available at  {$http$:$//cbcl.mit.edu/cbcl/software$-$datasets/FaceData2.html$}, consisting of 2429 gray-level images of faces (columns) with $19 \times 19$ pixels (rows), for which we set the factorization rank equal to $r=49$.}
\label{exNMF}
\end{center}
\end{figure*}

Nonnegative Matrix Factorization (NMF) consists in approximating a nonnegative matrix as the product of two low-rank nonnegative matrices~\cite{PT94, LS1}. More precisely, given an $m$-by-$n$ nonnegative matrix $M$ and a factorization rank $r$, we would like to find two nonnegative matrices $V$ and $W$ of dimensions $m$-by-$r$ and $r$-by-$n$ respectively such that 
\begin{displaymath}
M \approx V W.
\label{approx}
\end{displaymath}
This decomposition can be interpreted as follows: denoting by $M_{:j}$ the $j^{\textrm{th}}$ column of $M$, by $V_{:k}$ the $k^{\textrm{th}}$ column of $V$ and by $W_{kj}$ the entry of $W$ located at position $(k,j)$, we want
\begin{displaymath}
M_{:j} \; \approx \; \sum_{k=1}^r W_{kj} \; V_{:k}, 
\quad 1 \leq j \leq n,
\label{inter}
\end{displaymath}
so that each given (nonnegative) vector $M_{:j}$ is approximated by a nonnegative linear combination of $r$ nonnegative basis elements $V_{:k}$. Both the basis elements and the coefficients of the linear combinations have to be found. Nonnegativity of vectors $V_{:k}$ ensures that these basis elements belong to the same space $\mathbb{R}^m_+$ as the columns of $M$ and can then be interpreted in the same way. Moreover, the additive reconstruction due to nonnegativity of coefficients $W_{kj}$ leads to a \emph{part-based representation} \cite{LS1}: basis elements $V_{:k}$ will tend to represent common parts of the columns of $M$. 
For example, let each column of $M$ be a vectorized gray-level image of a face using (nonnegative) pixel intensities. The nonnegative matrix factorization of $M$ will generate a matrix $V$ whose columns are nonnegative basis elements of the original images, which can then be interpreted as images as well. Moreover, since each original face is reconstructed through a weighted sum of these basis elements, the latter will provide common parts extracted from the original faces, such as eyes, noses and lips. Figure~\ref{exNMF} illustrates this property of the NMF decomposition. 

One of the main challenges of NMF is to design fast and efficient algorithms generating the nonnegative factors. In fact, on the one hand, practitioners need to compute rapidly good factorizations for large-scale problems (e.g., in text mining or image processing); on the other hand, NMF is a NP-hard problem \cite{Vav} and we cannot expect to find a globally optimal solution in a reasonable computational time. 
This paper presents a general framework based on a multilevel strategy leading to faster initial convergence of NMF algorithms when dealing with data admitting a simple approximate low-dimensional representation (using linear transformations preserving nonnegativity), such as images. 
In fact, in these situations, a hierarchy of lower-dimensional problems can be constructed and used to compute efficiently approximate solutions of the original problem.  
Similar techniques have already been used for other dimensionality reduction tasks such as PCA~\cite{SFS09}. 

The paper is organized as follows: NMF is first formulated as an optimization problem and three well-known algorithms (ANLS, MU and HALS) are briefly presented. We then introduce the concept of multigrid/multilevel methods and show how and why it can be used to speed up NMF algorithms. Finally, we experimentally demonstrate the usefulness of the proposed technique on several standard image databases, and conclude with some remarks on limitations and possible extensions of this approach. \vspace{0.05cm}

\section{Algorithms for NMF} \label{algonmf}

NMF is typically formulated as a nonlinear optimization problem with an objective function measuring the quality of the low-rank approximation. In this paper, we consider the sum of squared errors:
\begin{equation}
\min_{\begin{array}{c} \, V \; \in \mathbb{R}^{m \times r}\\ W \in \mathbb{R}^{r \times n}\end{array}}   ||M-VW||_F^2  \quad \text{ s.t. } \quad  V \geq 0, \; W \geq 0, \label{NMF} \tag{NMF}
\end{equation}
i.e., use the squared Frobenius norm $||A||_F^2 = \sum_{i,j} A_{ij}^2$ of the approximation error. 
Since this standard formulation of \eqref{NMF} is NP-hard \cite{Vav}, most NMF algorithms focus on finding locally optimal solutions. In general, only convergence to stationary points of \eqref{NMF} (points satisfying the necessary first-order optimality conditions) is guaranteed.

\subsection{Alternating Nonnegative Least Squares (ANLS)} \label{anlsalgo}

Although \eqref{NMF} is a nonconvex problem, it is convex separately in each of the two factors $V$ and $W$, {\rm i.e.}, finding the optimal factor $V$ corresponding to a fixed factor $W$ reduces to a convex optimization problem, and vice-versa. More precisely, this convex problem corresponds to a nonnegative least squares (NNLS) problem, {\rm i.e.}, a least squares problem with nonnegativity constraints.  The so-called alternating nonnegative least squares (ANLS) algorithm for \eqref{NMF} minimizes (exactly) the cost function alternatively over factors $V$ and $W$ 
so that a stationary point of \eqref{NMF} is obtained in the limit \cite{Kim2}. 
\algsetup{indent=2em}
\begin{algorithm}[ht!]
\caption{Alternating Nonnegative Least Squares\label{ANLS}}
\begin{algorithmic}[1]
\REQUIRE Data matrix $M \in \mathbb{R}^{m \times n}_+$ and initial iterate $W \in \mathbb{R}^{r \times n}_+$.
\medskip
\WHILE {stopping criterion not met}
\STATE   $V \; \leftarrow \textrm{argmin}_{V \; \geq 0} ||M-VW||_F^2$;
\STATE   $W \leftarrow \textrm{argmin}_{W \geq 0} ||M-VW||_F^2$.
\ENDWHILE
\end{algorithmic}
\end{algorithm}
A frequent strategy to solve the NNLS subproblems is to use
active-set methods \cite{LH74} (see~\ref{AppA}) for which an efficient implementation\footnote{Available at \url{http://www.cc.gatech.edu/~hpark/}. Notice that an improved version based on a principal block pivoting method has been released recently, see \cite{KP08,KP11}, and for which our multilevel method is also applicable, see Section~\ref{ext}.} is described in \cite{VK04, Kim2}. 
We refer the reader to~\cite{Chen} for a survey about NNLS methods.

\subsection{Multiplicative Updates (MU)}

In \cite{LS2} Lee and Seung  propose multiplicative updates (MU) for \eqref{NMF} which guarantee nonincreasingness of the objective function (cf.\@ Algorithm~\ref{MU}). They also alternatively update $V$ for $W$ fixed and vice versa, using a technique which was originally proposed by Daube-Witherspoon and Muehllehner to solve nonnegative least squares problems \cite{DM86}. 
\algsetup{indent=2em}
\begin{algorithm}[ht!]
\caption{Multiplicative Updates\label{MU}}
\begin{algorithmic}[1]
\REQUIRE Data matrix $M \in \mathbb{R}^{m \times n}_+$ and initial iterates $(V,W) \in \mathbb{R}^{m \times r}_+ \times \mathbb{R}^{r \times n}_+$.
\medskip
\WHILE {stopping criterion not met} \vspace{0.1cm}
\STATE   $V \; \leftarrow V \; \circ \frac{[M W^T]}{[V (W W^T)]}$; 
\STATE   $W \leftarrow W \circ \frac{[V^T M]}{[(V^T V) W]}$.
\ENDWHILE
\end{algorithmic}
$\frac{[.]}{[.]}$ denotes the Hadamard (component-wise) division.
\end{algorithm}
The popularity of this algorithm came along with the popularity of NMF. Algorithm~\ref{MU} does not guarantee convergence to a stationary point 
(although it can be slightly modified in order to get this property~\cite{Lin2, Gil2}) and it has been observed to converge relatively slowly, see \cite{HHNP09} and the references therein.

\subsection{Hierarchical Alternating Least Squares (HALS)}

In ANLS, variables are partitioned at each iteration such that each subproblem is convex. However, the resolution of these convex NNLS subproblems is nontrivial and relatively expensive. If we optimize instead one single variable at a time, we get a simple univariate quadratic problem which admits a closed-form solution. 
Moreover, since the optimal value of each entry of $V$ (resp.\@ $W$) does not depend of the other entries of the same
column (resp.\@ row), 
one can optimize alternatively whole columns of $V$ and whole rows of $W$. 
This method was first proposed by Cichocki et al.\@ \cite{Cic, Cic4} and independently by \cite{diep, Gil1, LZ09}, and is herein referred to as Hierarchical Alternating Least Squares (HALS), see Algorithm~\ref{HALS}. 
\algsetup{indent=2em} 
\begin{algorithm}[ht!] 
\caption{Hierarchical Alternating Least Squares}\label{HALS} 
\begin{algorithmic}[1]
\REQUIRE Data $M \in \mathbb{R}^{m \times n}_+$ and initial iterates $(V,W) \in \mathbb{R}^{m \times r}_+ \times \mathbb{R}^{r \times n}_+$. 
\medskip
\WHILE {stopping criterion not met}
\STATE Compute $A = MW^T$ and $B = WW^T$.
 \FOR {$k = 1 \; : \; r$} \vspace{0.15cm}
			 \STATE $V_{:k} \leftarrow \max\Big(0, \frac{A_{:k} - \sum_{l = 1, l \neq k}^{r} V_{:l} B_{lk}}{B_{kk}}\Big)$; \vspace{0.1cm}
		   \ENDFOR			
\STATE Compute $C = V^TM$ and $D = V^TV$.
 \FOR {$k = 1 \; : \; r$} \vspace{0.15cm}
			 \STATE $W_{k:} \leftarrow \max\Big(0, \frac{C_{k:} - \sum_{l = 1, l \neq k}^{r} D_{kl} W_{l:}}{D_{kk}}\Big)$; \vspace{0.1cm}
		   \ENDFOR	
\ENDWHILE
\end{algorithmic}
\end{algorithm}
Under some mild assumptions, every limit point is a stationary point of \eqref{NMF}, see \cite{Gillis}. 

\section{Multigrid Methods} \label{multigrid}

In this section, we briefly introduce multigrid methods. 
The aim is to give the reader some insight on these techniques in order to comprehend their applications for NMF. We refer the reader to~\cite{Bra, Bran, Bri, Trot} and the references therein for detailed discussions on the subject.\\

Multigrid methods were initially used to develop fast numerical solvers for boundary value problems. 
Given a differential equation on a continuous domain with boundary conditions, the aim is to find an approximation of a \emph{smooth} function $f$ satisfying the constraints. In general, the first step is to discretize the continuous domain, {\rm i.e.}, choose a set of points (a \emph{grid}) where the function values will be approximated. Then, a numerical method (\rm{e.g.}, finite differences, finite elements)  translates the continuous problem into a (square) system of linear equations:
\begin{equation} \label{lins}
\textrm{find } x \in \mathbb{R}^{n} \quad 
\textrm{ s.t. } \quad     Ax = b,\quad \textrm{ with } A \in \mathbb{R}^{n \times n},\; b \in \mathbb{R}^{n},
\end{equation}
where the vector $x$ will contain the approximate values of $f$ on the grid points. 
Linear system \eqref{lins} can be solved either by direct methods (\rm{e.g.}, Gaussian elimination) or iterative methods (\rm{e.g.}, Jacobi and Gauss-Seidel iterations). Of course, the computational cost of these methods depends on the number of points in the grid, which leads to a trade-off between precision (number of points used for the discretization) and computational cost.\\

Iterative methods update the solution at each step 
and hopefully converge to a solution of \eqref{lins}. 
Here comes the utility of multigrid: instead of working on a fine grid during all iterations, the solution is  \emph{restricted} to a coarser grid\footnote{Standard multigrid techniques actually restrict the residual instead of the solution, see the discussion in Section~\ref{conv}.} on which the iterations are cheaper. Moreover, the smoothness of function $f$ allows to recover its low-frequency components faster on coarser grids.  
Solutions of the coarse grid are then \emph{prolongated} to the finer grid and iterations can continue (higher frequency components of the error are reduced faster). Because the initial guess generated on the coarser grid is a good approximation of the final solution, less iterations are needed to converge on the fine (expensive) grid.
Essentially, multigrid methods make iterative methods more efficient, {\rm i.e.}, accurate solutions are obtained faster.\\

More recently, these same ideas have been applied to a broader class of problems, \rm{e.g.}, multiscale optimization with trust-region methods~\cite{Sart} and multiresolution techniques in image processing~\cite{Ter}.

\section{Multilevel Approach for NMF} \label{multiNMF}

The three algorithms presented in Section~\ref{algonmf} (ANLS, MU and HALS) are iteratively trying to find a stationary point of \eqref{NMF}. Indeed, most practical NMF algorithms are \emph{iterative methods}, such as projected gradient methods~\cite{Lin} and Newton-like methods~\cite{Zdu1, Dhi2} (see also~\cite{Ber, CAZP09, Cic2, diep} and the references therein).
In order to embed these algorithms in a multilevel strategy, one has to define the different levels and describe how variables and data are transferred between them. 
In this section, we first present a general description of the multilevel approach for NMF algorithms, and then apply it to image datasets.

\subsection{Description}

Let each column of the matrix $M$ be a element of the dataset (e.g., a vectorized image) belonging to $\mathbb{R}^m_+$. 
We define the restriction operator $\mathcal{R}$ as a linear operator
\[
\mathcal{R} : \mathbb{R}^m_+ \rightarrow \mathbb{R}^{m'}_+ : x \rightarrow \mathcal{R}(x) = R x,
\]
with $R \in \mathbb{R}^{m' \times m}_+$ and $m' < m$, 
and the prolongation $\mathcal{P}$ as a linear operator
\[
\mathcal{P} : \mathbb{R}^{m'}_+ \rightarrow \mathbb{R}^{m}_+ : y \rightarrow \mathcal{P}(y) = P y,
\]
with $P \in \mathbb{R}^{m \times m'}_+$. Nonnegativity of matrices $R$ and $P$ is a sufficient condition to preserve nonnegativity of the solutions when they are transferred from one level to another. In fact, in order to generate nonnegative solutions, one requires
\[
\mathcal{R}(x) \geq 0, \, \forall x \geq 0 \; \textrm{ and } \; \mathcal{P}(y) \geq 0, \, \forall y \geq 0.
\]
We also define the corresponding transfer operators on matrices, operating columnwise:
\[
\mathcal{R}([x_1 \, x_2 \dots x_n]) = [\mathcal{R}(x_1) \, \mathcal{R}(x_2) \dots \mathcal{R}(x_n)], \textrm{ and}
\]
\[
\mathcal{P}([y_1 \, y_2 \dots y_n]) = [\mathcal{P}(y_1) \, \mathcal{P}(y_2) \dots \mathcal{P}(y_n)],
\]
for $x_i \in \mathbb{R}^{m}_+, y_i \in \mathbb{R}^{m'}_+, 1 \leq i \leq n$.\\

In order for the multilevel strategy to work, information lost when transferring from one level to another must be limited, i.e., data matrix $M$ has to be well represented by $\mathcal{R}(M)$ in the lower dimensional space, which means that the reconstruction $\mathcal{P}(\mathcal{R}(M))$ must be close to $M$. From now on, we say that $M$ is smooth with respect to $\mathcal{R}$ and $\mathcal{P}$ if and only if
\[
s_M = \frac{|| M - \mathcal{P}(\mathcal{R}(M))||_F}{||M||_F} \quad \textrm{ is small }.
\]
Quantity $s_M$ measures how well $M$ can be mapped by $\mathcal{R}$ into a lower-dimensional space, then brought back by $\mathcal{P}$, and still be a fairly good approximation of itself. \\

Based on these definitions, elaborating a multilevel approach for NMF is straightforward:\\ 

\begin{enumerate}

\item We are given $M \in \mathbb{R}^{m \times n}_+$ and $(V_0,W_0) \in \mathbb{R}^{m \times r}_+ \times \mathbb{R}^{r \times n}_+$;

\item Compute $M' = \mathcal{R}(M) = RM \in \mathbb{R}^{m' \times n}_+$ and $V_0' = \mathcal{R}(V_0) = RV_0 \in \mathbb{R}^{m' \times r}_+$, i.e., restrict the elements of your dataset and the basis elements of the current solution to a lower dimensional space;

\item Compute a rank-$r$ NMF $(V',W)$ of $M'$ using $(V_0',W_0)$ as initial matrices, i.e., 
\[
V'W \approx M' = \mathcal{R}(M).
\]
This can be done using any NMF iterative algorithm or, even better, using the multilevel strategy recursively (cf.\@ Section~\ref{cycle}).

\item Since 
\[
M \approx \mathcal{P}(\mathcal{R}(M)) = \mathcal{P}(M') \approx \mathcal{P}(V'W) = PV'W = \mathcal{P}(V') W = VW,
\]
where $V$ is the prolongation  of $V'$,  $(V,W)$ is a good initial estimate for a rank-$r$ NMF of $M$, \emph{provided that} $M$ is smooth with respect to $\mathcal{R}$ and $\mathcal{P}$ (i.e., $s_M$ is small) and that $V'W$ is a good approximation of $M' = \mathcal{R}(M)$ (i.e., $||M'-V'W||_F$ is small); in fact, 
\begin{eqnarray*}
||M - \mathcal{P}(V')W||_F & \leq & ||M - \mathcal{P}(\mathcal{R}(M))||_F + ||\mathcal{P}(\mathcal{R}(M)) - \mathcal{P}(V'W)||_F \\
& \leq & s_M ||M||_F + ||\mathcal{P}(\mathcal{R}(M)-V'W)||_F \\
& \leq & s_M ||M||_F + ||P||_F \, ||\mathcal{R}(M)-V'W||_F.
\end{eqnarray*}
\item Further improve the solution $(V,W)$ using any NMF iterative algorithm.\\

\end{enumerate}

Computations needed at step~3 are cheap (since $m' < m$) and, moreover, the low-frequency components of the error\footnote{The low-frequency components refer to the parts of the data which are well-represented on coarse levels.} are reduced faster on coarse levels (cf.\@ Section~\ref{smpro}). Therefore this strategy is expected to accelerate the convergence of NMF algorithms.

We now illustrate this technique on image datasets, more precisely, on two-dimensional gray-level images. In general, images are composed of several smooth components, i.e., regions where pixel values are similar and change continuously with respect to their location (e.g., skin on a face or, the pupil or sclera\footnote{The white part of the eye.} of an eye). In other words, a pixel value can often be approximated using the pixel values of its neighbors. This observation can be used to define the transfer operators (Section~\ref{operators}). For the computation of a NMF solution, the multilevel approach can be used recursively; three strategies (called multigrid cycles) are described in Section~\ref{cycle}. Finally, numerical results are reported in Section~\ref{Res}. 

\subsection{Coarse Grid and Transfer Operators} \label{operators}

\begin{figure*}[ht!]
\begin{center}
\includegraphics[width=12cm]{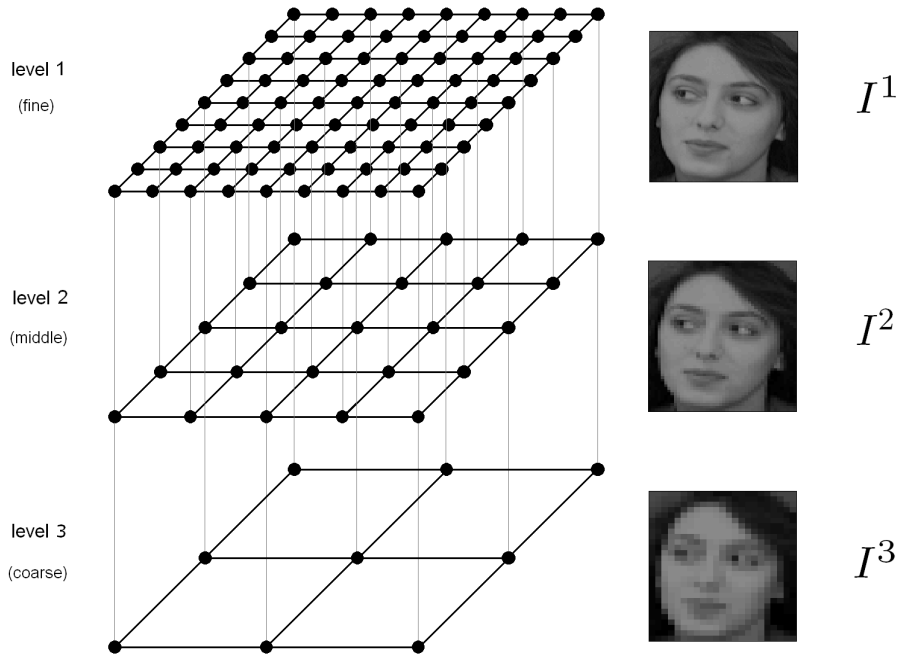}
\caption{Multigrid Hierarchy. Schematic view of a grid definition for image processing (image from ORL face database, cf.\@ Section~\ref{Res}).}
\label{multires}
\end{center}
\end{figure*}

A crucial step of multilevel methods is to define the different levels and the transformations (operators) between them. 
Figure \ref{multires} is an illustration of a standard \emph{coarse grid} definition: 
we note $I^1$ the matrix of dimension $(2^{a}+1) \times (2^{b}+1)$ representing the initial image 
and $I^l$ the matrix of dimension $(2^{a-l+1}+1) \times (2^{b-l+1}+1)$ representing the image at level $l$ obtained by keeping, in each direction, only one out of every two points of the grid at the preceding level, {\rm i.e.}, $I^{l-1}$.

\begin{figure}[ht!]
\begin{center}
\includegraphics[width=7cm]{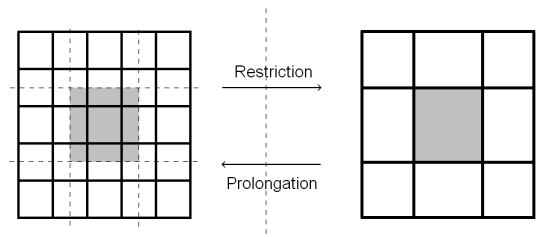}
\caption{Restriction and Prolongation.}
\label{respro}
\end{center}
\end{figure}

The transfer operators describe how to transform the images when going from finer to coarser levels, and vice versa, {\rm i.e.}, how to compute the values (pixel intensities) of the image $I^l$ using values from image $I^{l-1}$ at the finer level (restriction) or from image $I^{l+1}$ at the coarser level (prolongation). For the \emph{restriction}, the \emph{full-weighting} operator is a standard choice: values of the coarse grid points are the weighted average of the values of their neighbors on the fine grid (see Figure~\ref{respro} for an illustration). 
\noindent Noting $I_{i,j}^l$ the intensity of the pixel $(i,j)$ of image $I^l$, 
it is defined as follows: 
\begin{align}
I_{i,j}^{l+1} =  & \frac{1}{16} \Big[ I_{2i-1,2j-1}^l + I_{2i-1,2j+1}^l + I_{2i+1,2j-1}^l + I_{2i+1,2j+1}^l \nonumber \\
& + 2 (I_{2i,2j-1}^l + I_{2i-1,2j}^l + I_{2i+1,2j}^l + I_{2i,2j+1}^l) + 4 I_{2i,2j}^l \Big],\label{restr} 
\end{align}
except on the boundaries of the image (when $i=0$, $j = 0$, $i= 2^{a-l+1}$ and/or $j= 2^{b-l+1}$) where the weights are adapted correspondingly. For example, to restrict a $3 \times 3$ image to a $2 \times 2$ image, $\mathcal{R}$ is defined with
\begin{displaymath}
R = \frac{1}{9} 
\left( \begin{array}{cccccccccc}
4 & 2 & 0 & 2 & 1 &0 &0 &0 &0 \\
0 & 2 & 4 & 0 & 1 &2 &0 &0 &0 \\
0 & 0 & 0 & 2 & 1 &0 &4 &2 &0 \\
0 & 0 & 0 & 0 & 1 &2 &0 &2 & 4\end{array} \right),
\end{displaymath}
($3\times 3$ images needing first to be vectorized to vectors in $\mathbb{R}^9$, by concatenation of either columns or rows).

For the \emph{prolongation}, we set the values on the fine grid points as the average of the values of their neighbors on the coarse grid:
\begin{equation} \label{prol}
I_{i,j}^{l} = {\text{mean}}_{\begin{subarray}{l} i' \in \textrm{rd}(i/2) \\ j' \in \textrm{rd}(j/2) \end{subarray}}  \, \Big(I_{i',j'}^{l+1}\Big),
\end{equation}
where
\[ \textrm{rd}(k/2) = \left\{ 
\begin{array}{ll}
\{k/2\}																					& \textrm{$k$ even,}\\
\{(k-1)/2, (k+1)/2 \}	 & \textrm{$k$ odd}.
\end{array} \right.
\]
For example, to prolongate a $2 \times 2$ image to a $3 \times 3$ image, $\mathcal{P}$ is defined with
\begin{displaymath}
P^T = \frac{1}{4} 
\left( \begin{array}{cccccccccc}
     4   &  2    & 0 &    2   &  1 &    0  &   0 &    0 &    0 \\
     0   &  2   &  4 &    0   &  1 &    2  &   0 &    0 &    0\\
     0   &  0  &   0 &    2   &  1 &    0  &   4 &    2 &    0\\
     0   &  0 &    0 &    0 &    1 &    2 &    0 &    2 &    4 \end{array} \right).
\end{displaymath}

Note that these transformations clearly preserve nonnegativity.

\subsection{Multigrid Cycle} \label{cycle}

Now that grids and transfer operators are defined, we need to choose the procedure that is applied at each grid level as it moves through the grid hierarchy. In this section, we propose three different approaches: nested iteration, V-cycle and full multigrid cycle. \\

In our setting, the transfer operators only change the number of rows $m$ of the input matrix $M$, {\rm i.e.}, the number of pixels in the images of the database: the size of the images  are approximatively four times smaller between each level: $m' \approx \frac{1}{4}m$. When the number of images in the input matrix is not too large, i.e.\@, when $n \ll m$, the computational complexity per iteration of the three algorithms (ANLS, MU and HALS) is close to being proportional to $m$ (cf.~\ref{AppA}), and the iterations will then be approximately four times cheaper (see also Section~\ref{sizeM}). 
A possible way to allocate the time spent at each level is to allow the same number of iterations at each level, which seems to give good results in practice. Table~\ref{table1} shows the time spent and the corresponding number of iterations performed at each level. 
\begin{table}[ht!]
\caption{Number of iterations performed and time spent at each level when allocating among $L$ levels a total computational budget $T$, corresponding to $4k$ iterations at the finest level.} 
\begin{center}
\renewcommand\arraystretch{1.2}
\begin{tabular}{c|c|c|c|c|c|c}
	&   Level $1$   &  Level $2$ & \dots 	&  Level $L-1$ &  Level $L$ & Total  \\
	&   (finer)    &   & \dots 	&   &  (coarser)   & \\
\hline  
\vspace{0.1cm}  \textit{$\sim$  \# iterations}& $3k$ & $3k$   & \dots 	& $3k$ & $4k$  & $(3L+1)k$ \\ 
 \textit{time} & $\frac{3}{4}T$   & $\frac{3}{16} T$ & \dots & $\frac{3}{4^{\textrm{L-1}}} T$	& $\frac{1}{4^{\textrm{L-1}}} T$ & $T$
\end{tabular}
\end{center}
\label{table1}
\end{table}

Note that the transfer operators require $\mathcal{O}(mn)$ operations and, since they are only performed once between each level, their computational cost can be neglected (at least for $r \gg 1$ and/or when a sizable amount of iterations are performed).

\subsubsection{Nested Iteration (NI)}

\begin{figure*}[ht!]
\begin{center}
\includegraphics{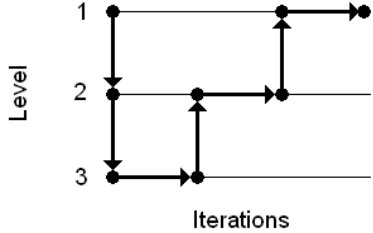} 
\caption{Nested Iteration. Transition between different levels for nested iteration.}
\label{ni}
\end{center}
\end{figure*}

To initialize NMF algorithms, we propose to factorize the image at the coarsest resolution and then use the solution as initial guess for the next (finer) resolution. This is referred to as \emph{nested iteration}, see Figure~\ref{ni} for an illustration with three levels and Algorithm~\ref{nit} for the implementation. The idea is to start off the final iterations at the finer level with a better initial estimate, thus reducing the computational time required for the convergence of the iterative methods on the fine grid. The number of iterations and time spent at each level is chosen according to Table~\ref{table1}, i.e., three quarters of the alloted time for iterations at the current level preceded by one quarter of the time for the recursive call to the immediately coarser level.
\algsetup{indent=2em}
\begin{algorithm}[ht!]
\caption{Nested Iteration}\label{nit}
\begin{algorithmic}[1]
\REQUIRE $L \in \mathbb{N}$ (number of levels), $M \in \mathbb{R}^{m \times n}_+$ (data matrix), 
$(V_0,W_0) \in \mathbb{R}^{m \times r}_+ \times \mathbb{R}^{r \times n}_+$ (initial matrices) and $T \geq 0$ (total time allocated to the algorithm).
\ENSURE $(V,W) \geq {0}$ s.t. $VW \approx M$.
\medskip
\IF {$L = 1$}
	 \STATE $[V,W] = \textrm{NMF algorithm}(M,V_0,W_0,T)$;
\ELSE
   \STATE $M' = \mathcal{R}(M)$; $V_0' = \mathcal{R}(V_0)$;
	 \STATE $[V',W] = \textrm{Nested Iteration}(L-1,M',V_0',W_0,T/4)$;
	 \STATE $V = \mathcal{P}(V')$;  \vspace{0.1cm}
	 \STATE $[V,W] = \textrm{NMF algorithm}(M,V,W,3T/4)$;
\ENDIF			
\end{algorithmic}
\end{algorithm}

\begin{remark}
When the ANLS algorithm is used, the prolongation of $V'$ does not need to be computed since that algorithm only needs an initial value for the $W$ iterate. Note that this can be used in principle to avoid computing any prolongation, by setting $V$ directly as the optimal solution of the corresponding NNLS problem.
\end{remark}

\subsubsection{V--Cycle (VC)}

\begin{figure*}[ht!]
\begin{center}
\includegraphics{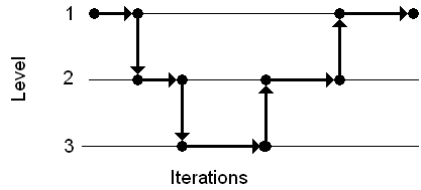} 
\caption{V-cycle. Transition between different levels for V-cycle.}
\label{vcyc}
\end{center}
\end{figure*}

One can often empirically observe that multilevel methods perform better if a few iterations are performed at the fine level immediately before going to coarser levels. This is partially explained by the fact that these first few iterations typically lead to a relatively important decrease of the objective function, at least compared to subsequent iterations. A simple application of this strategy is referred to as V-cycle and is illustrated on Figure~\ref{vcyc} with three levels; see Algorithm~\ref{vcy} for the implementation. Time allocation is as follows: one quarter of the alloted time is devoted to iterations at the current level, followed by one quarter of the time for the recursive call to the immediately coarser level, and finally one half of the time again for iterations at the current level (we have therefore three quarters of the total time spent for iterations at current level, as for nested iteration).

\algsetup{indent=2em}
\begin{algorithm}[ht!]
\caption{V-cycle}\label{vcy}
\begin{algorithmic}[1]
\REQUIRE $L \in \mathbb{N}$ (number of levels), $M \in \mathbb{R}^{m \times n}_+$ (data matrix), 
$(V_0,W_0) \in \mathbb{R}^{m \times r}_+ \times \mathbb{R}^{r \times n}_+$ (initial matrices) and $T \geq 0$ (total time allocated to the algorithm).
\ENSURE $(V,W) \geq {0}$ s.t. $VW \approx M$.
\medskip
\IF {$L = 1$}
	 \STATE $[V,W] = \textrm{NMF algorithm}(M,V_0,W_0,T)$;
\ELSE
	 \STATE $[V,W] = \textrm{NMF algorithm}(M,V_0,W_0,T/4)$;
   \STATE $M' = \mathcal{R}(M)$; $V'  = \mathcal{R}(V)$;
	 \STATE $[V',W] = \textrm{V-cycle}(L-1,M',V',W,T/4)$;
	 \STATE $V = \mathcal{P}(V')$; \vspace{0.1cm}
	 \STATE $[V,W] = \textrm{NMF algorithm}(M,V,W,T/2)$;
\ENDIF			
\end{algorithmic}
\end{algorithm}

\subsubsection{Full Multigrid (FMG)}

Combining ideas of nested iteration and V-cycle leads to a full multigrid cycle defined recursively as follows: at each level, a V-cycle is initialized with 
the solution obtained at the underlying level using a full-multigrid cycle. 
This is typically the most efficient multigrid strategy \cite{Trot}. In this case, we propose to partition the time as follows ($T$ is the total time): $\frac{T}{4}$ for the initialization (call of the full multigrid on the underlying level) and $\frac{3T}{4}$ for the V-cycle at the current level (cf.\@ Algorithm~\ref{fmg}).
\algsetup{indent=2em}
\begin{algorithm}[ht!]
\caption{Full Multigrid}\label{fmg}
\begin{algorithmic}[1]
\REQUIRE $L \in \mathbb{N}$ (number of levels), $M \in \mathbb{R}^{m \times n}_+$ (data matrix), 
$(V_0,W_0) \in \mathbb{R}^{m \times r}_+ \times \mathbb{R}^{r \times n}_+$ (initial matrices) and $T \geq 0$ (total time allocated to the algorithm).
\ENSURE $(V,W) \geq {0}$ s.t. $VW \approx M$.
\medskip
\medskip
\IF {$L = 1$}
	 \STATE $[V,W] = \textrm{NMF algorithm}(M,V_0,W_0,T)$;
\ELSE
	\STATE $V' = \mathcal{R}(V_0)$; $M' = \mathcal{R}(M)$; * 
	\STATE $[V',W] = \textrm{Full Multigrid}(L-1,M',V',W_0,T/4)$;
	\STATE $V = \textrm{prolongation}(V')$;
  \STATE $[V,W] = \textrm{V-cycle}(L,M,V,W,3T/4)$; 
\ENDIF \\
\medskip
\footnotesize*Note that the restrictions of $M$ should be computed only once for each level and saved as global variables so that the call of the V-cycle (step~7) does not have to recompute them.
\end{algorithmic}
\end{algorithm}

\subsection{Smoothing Properties} \label{smpro}

We explained why the multilevel strategy was potentially able to accelerate iterative algorithms for NMF: cheaper computations and smoothing of the error on coarse levels.
Before giving extensive numerical results in Section~\ref{Res}, we illustrate these features of multilevel methods on the ORL face database.

\begin{figure*}[ht!]
\begin{center}
\includegraphics[width=\textwidth]{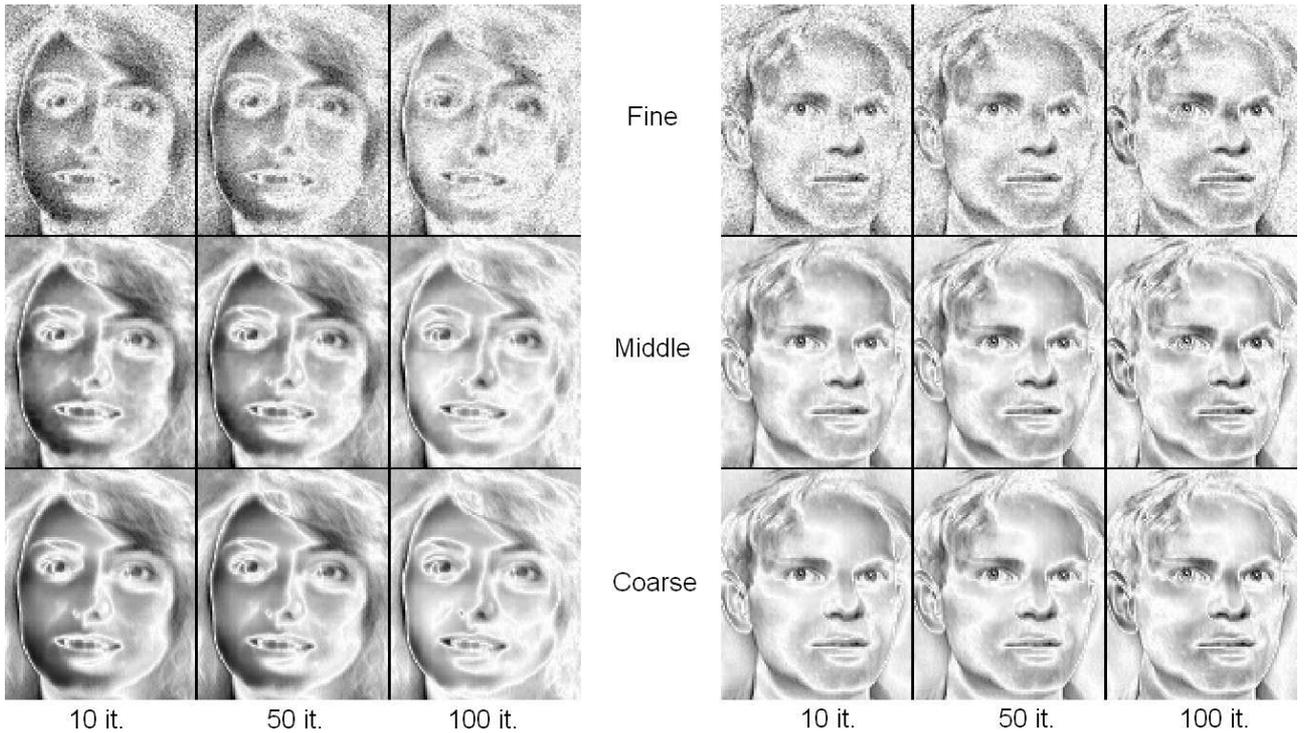}
\caption{Smoothing on Coarse Levels. Example of the smoothing properties of the multilevel approach on the ORL face database. 
Each image represents the absolute value of the approximation error (black tones indicate a high error) of one of  two faces from the ORL face database. 
These approximations are the prolongations (to the fine level) of the solutions obtained using the multiplicative updates on a single level, with factorization rank $r=40$ and the same initial matrices.  
From top to bottom: level 1 (fine), level 2 (middle) and level 3 (coarse); from left to right: 10 iterations, 50 iterations and 100 iterations.}
\label{smootherror}
\end{center}
\end{figure*}

\begin{figure*}[ht!]
\begin{center}
\includegraphics[width=\textwidth]{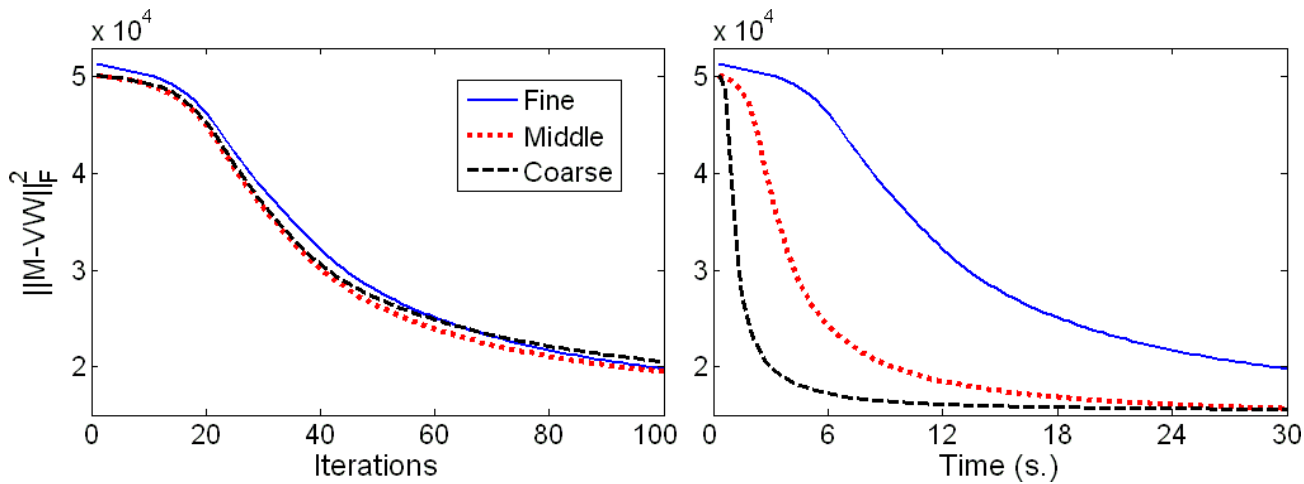}
\caption{Evolution of the error on each level, after prolongation on the fine level,  with respect to (left) the number of iterations performed and (right) the computational time. Same setting as in Figure~\ref{smootherror}.} 
\label{errsmo}
\end{center}
\end{figure*}

Comparing three levels, 
Figure~\ref{smootherror} displays the error (after prolongation to the fine level) for two faces and for different numbers of iterations (10, 50 and 100) using MU.   
Comparing the first row and the last row of Figure~\ref{smootherror}, it is clear that, in this example, the multilevel approach allows a significant smoothing of the error. After only 10 iterations, the error obtained with the prolongated solution of the coarse level is already smoother and smaller (see Figure~\ref{errsmo}), while it is computed much faster.  

Figure~\ref{errsmo} gives the evolution of the error with respect to the number of iterations performed (left) and with respect to computational time (right). In this example, the initial convergence on the three levels is comparable, while the computational cost is much cheaper on coarse levels. In fact, compared to the fine level, the middle (resp.\@ coarse) level is approximately 4 (resp.\@ 16) times cheaper. 

\section{Computational Results} \label{Res}

To evaluate the performances of our multilevel approach, we present some numerical results for several standard image databases described in Table~\ref{dipr}. 

\begin{table}[ht!]
\caption{Image datasets.}
\begin{center}
\begin{tabular}{|c|c|c|c|c|}
\hline
Data &             $\#$ pixels &  m  &  n  & r\\ \hline \hline
ORL face$^1$    &   $112 \times 92$  & 10304 & 400 & 40 \\ 
Umist face$^2$  &   $112 \times 92$  & 10304 & 575 & 40 \\ 
Iris$^3$          &   $960 \times 1280$ & 1228800 & 8 & 4 \\ 
Hubble Telescope \cite{hubble} &  $128 \times 128$  & 16384  & 100  & 8 \\ 
\hline
\end{tabular} \\
\begin{flushleft}
\footnotesize$^1$
\url{http://www.cl.cam.ac.uk/research/dtg/attarchive/facedatabase.html}\\
$^2$ \url{http://www.cs.toronto.edu/~roweis/data.html}\\ 
$^3$ \url{http://www.bath.ac.uk/elec-eng/research/sipg} 
\end{flushleft}
\end{center}
\label{dipr}
\end{table}
For each database, the three multigrid cycles (NI, V-cycle and FMG) of our multilevel strategy are tested using 100 runs initialized with the same  random matrices for the three algorithms (ANLS, MU and HALS), with a time limit of 10 seconds. 
All algorithms have been implemented in \matlab \, 7.1 (R14) and tested on a 3~GHz Intel$^{\textrm{\textregistered}}$ Core{\texttrademark}2 Dual CPU PC.

\subsection{Results}

Tables~\ref{errortime1}, \ref{errortime2} and \ref{errortime3} give the mean error attained within 10 seconds using the different approaches. 
\begin{center}
\begin{table}[ht!] 
\caption{Comparison of the mean error on the 100 runs with ANLS.}
\begin{center}
\begin{tabular}{|c|c|c|c|c|c|}
\hline
    & \# lvl 	 &  ORL & Umist & Iris & Hubble   \\ \hline
NMF & 1 	 &  14960  &  26013 & 28934   &  24.35            \\ 
\hline
NI & 2 	 &  14683   &  25060 &  27834   &  \textbf{15.94}  \\ 
	  & 3  &  14591   &  24887 &  27572   &  16.93    \\
	  & 4  &  14580   &  24923 &  27453   &  17.20   \\
	  \hline
VC  & 2  & 14696    & 25195  &  27957   &  16.00   \\ 
	  & 3  & 14610    & 24848  &  27620   &  16.12  \\ 
	  & 4  & 14599    & 24962  &  27490   &  16.10  \\ 
\hline
FMG & 2  &  14683   & 25060  &  27821   &  16.10      \\ 
	  & 3  &	14516   & 24672  &  27500   &  16.56      \\ 
	  & 4  &	\textbf{14460}   & \textbf{24393}  &  \textbf{27359}   &  16.70       \\ 
\hline
\end{tabular}
\label{errortime1}
\end{center}
\end{table}
\end{center}
\begin{center}
\begin{table}[ht!] 
\caption{Comparison of the mean error on the 100 runs with MU.}
\begin{center}
\begin{tabular}{|c|c|c|c|c|c|}
\hline
    & \# lvl 	 &  ORL & Umist & Iris & Hubble   \\ \hline
NMF & 1 	 & 34733  &  131087 & 64046   &  21.68   \\ 
\hline
NI & 2 	 &  23422   &  87966 & 37604    & 22.80   \\ 
	  & 3  &  20502   &  67131 & 33114    & 18.49    \\
	  & 4  &  19507   &  59879 & 31146    & 16.19  \\
	  \hline
VC  & 2  &  23490   &  90064 &  36545   & 10.62    \\ 
	  & 3  &  20678   &  69208 &  32086   & 9.77   \\ 
	  & 4  &  19804   &  62420 &  30415   & \textbf{9.36}   \\ 
\hline
FMG & 2  &  23422   &  87966 &  37504   & 22.91    \\ 
	  & 3  &	19170   &  58469 &  32120   & 15.06     \\ 
	  & 4  &	\textbf{17635}   &  \textbf{46570} &  \textbf{29659}   & 11.71     \\ 
\hline
\end{tabular}
\label{errortime2}
\end{center}
\end{table}
\end{center}
\begin{center}
\begin{table}[ht!] 
\caption{Comparison of the mean error on the 100 runs with HALS.}
\begin{center}
\begin{tabular}{|c|c|c|c|c|c|}
\hline
    & \# lvl 	 &  ORL & Umist & Iris & Hubble   \\ \hline
NMF & 1 	 & 15096  &  27544 &  31571   &  17.97            \\ 
\hline
NI & 2 	 &  14517   &  25153 &  29032   &  17.37    \\ 
	  & 3  &  14310   &  24427 &  28131   &  16.91    \\
	  & 4  &  14280   &  24256 &  27744   &  16.92    \\
	  \hline 
VC  & 2  &  14523   &  25123 &  28732   &  17.37    \\ 
	  & 3  &  14339   &  24459 &  28001   &  17.02    \\ 
	  & 4  &  14327   &  24364 &  27670   &  17.04    \\ 
\hline
FMG & 2  &  14518   &  25153 &  29120   &  17.39    \\ 
	  & 3  &	14204   &  23950 &  27933   &  \textbf{16.69}    \\ 
	  & 4  &	\textbf{14107}   &  \textbf{23533} &  \textbf{27538}   &  16.89    \\ 
\hline
\end{tabular}
\label{errortime3}
\end{center}
\end{table}
\end{center}
In all cases, the multilevel approach generates much better solutions than the original NMF algorithms, indicating that it is able to accelerate their convergence. The full multigrid cycle is, as expected, the best strategy while nested iteration and V-cycle give comparable performances. 
We also observe that the additional speedup of the convergence when the number of levels is increased from 3 to 4 is less significant; it has even a slightly negative effect in some cases. In general, the `optimal' number of levels will depend on the smoothness and the size of the data, and on the algorithm used (cf.\@ Section~\ref{sizeM}).\\

\begin{figure*}[ht!]
\begin{center}
\includegraphics[width=\textwidth]{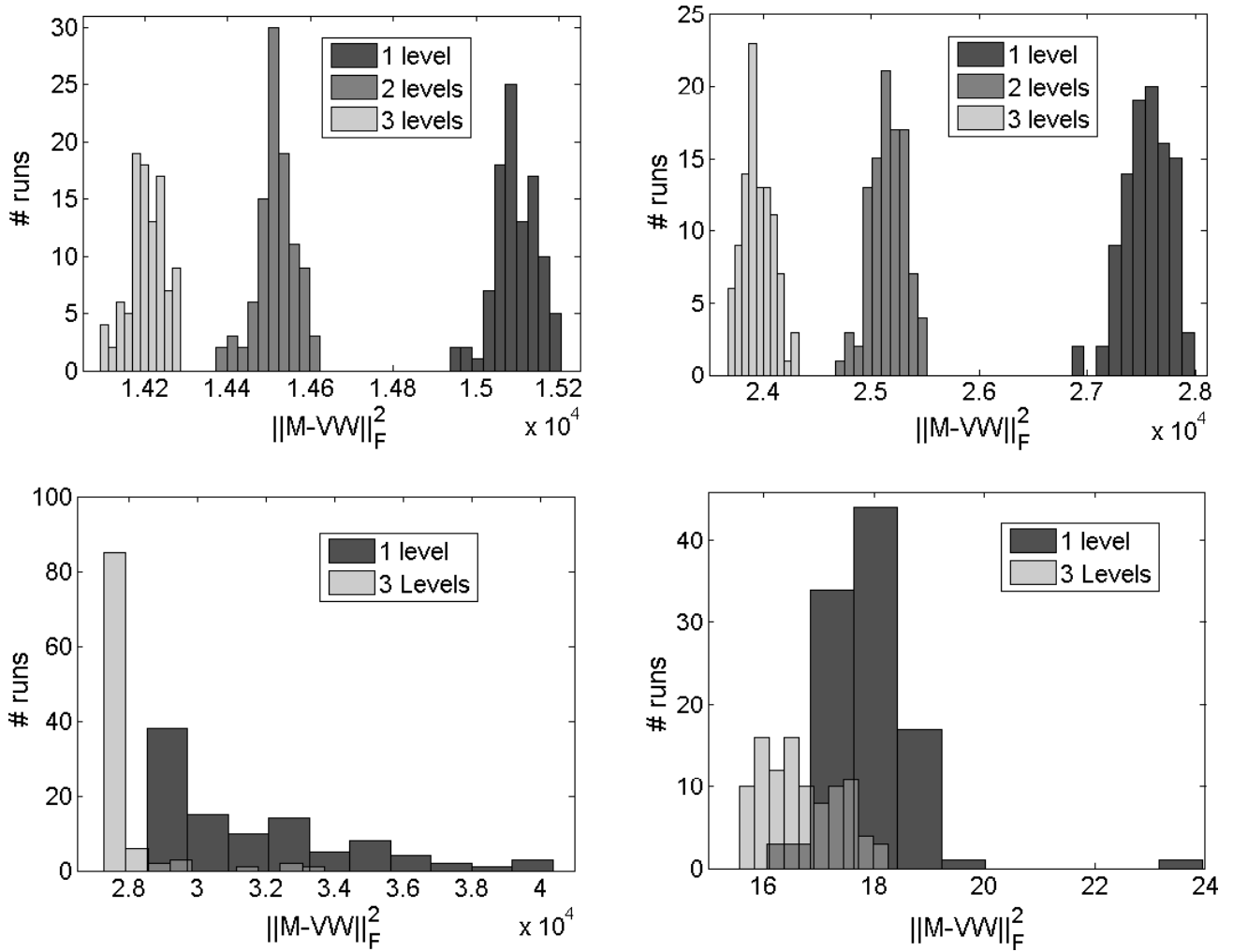}
\caption{Distribution of the error among the 100 random initializations using the HALS algorithm with a full multigrid cycle: (top left) ORL, (top right) Umist, (bottom left) Iris, and (bottom right) Hubble.} 
\label{histo}
\end{center}
\end{figure*}

HALS combined with the the full multigrid cycle is one of the best strategies. 
Figure~\ref{histo} displays the distribution of the errors for the different databases in this particular case. 
For the ORL and Umist databases, the multilevel strategy is extremely efficient: all the solutions generated with 2 and 3 levels are better than the original NMF algorithm. For the Iris and Hubble databases, the difference is not as clear.  
The reason is that the corresponding NMF problems are `easier' because the factorization rank $r$ is smaller. Hence the algorithms converge faster to stationary points, and the distribution of the final errors is more concentrated. 

\begin{figure*}[ht!]
\begin{center}
\includegraphics[width=\textwidth]{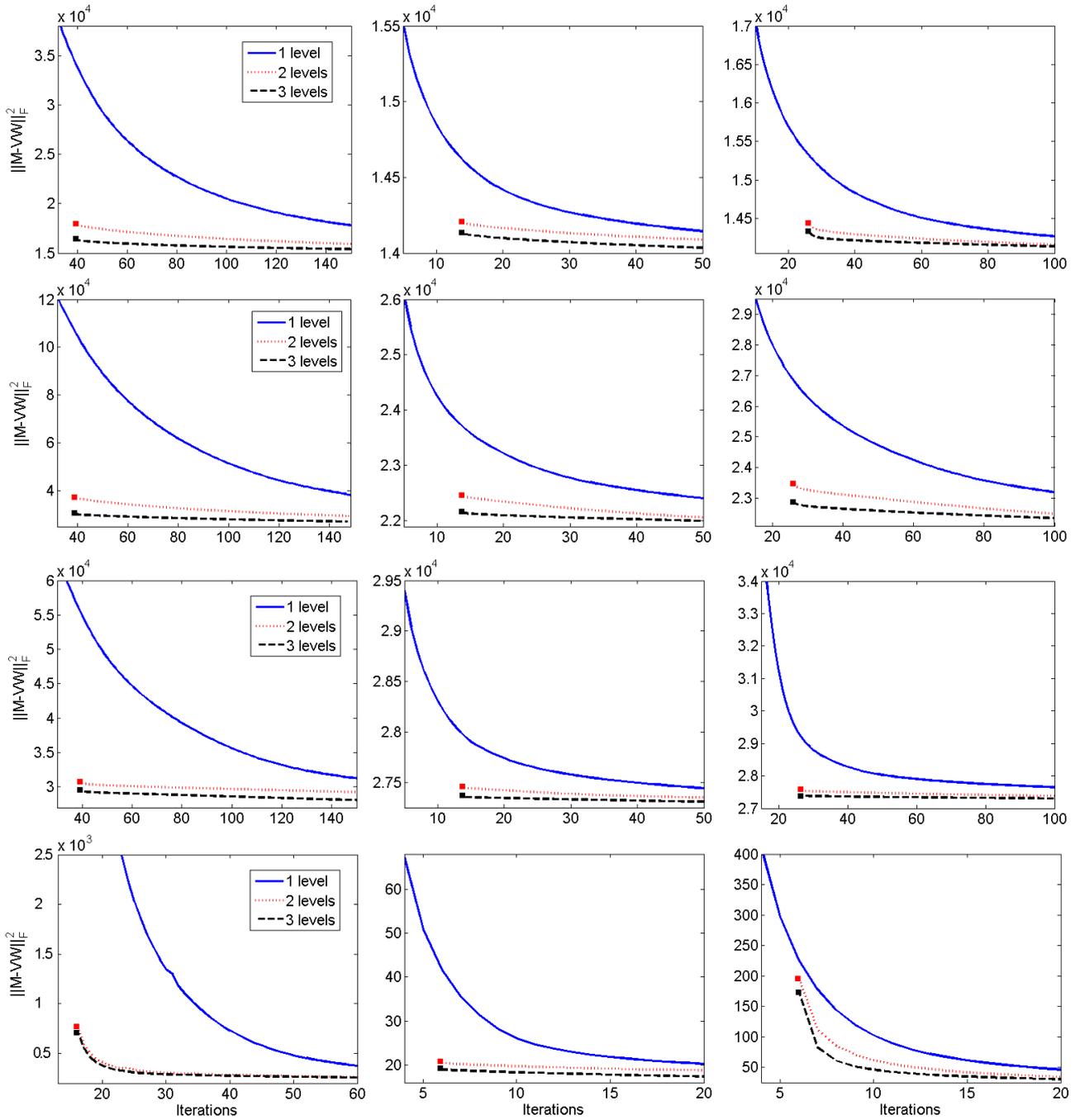}
\caption{Evolution of the objective function. From left to right~: MU, ANLS and HALS. From top to bottom: ORL, Umist, Iris and Hubble databases. \emph{1 level} stands for the standard NMF algorithms. The initial points for the curves \emph{2 levels} and \emph{3 levels} are the prolongated solutions obtained on the coarser levels using nested iteration, cf.~Section~\ref{cycle}. All algorithms were initialized with the same random matrices.}
\label{plotAll}
\end{center}
\end{figure*}

In order to visualize the evolution of the error through the iterations, Figure~\ref{plotAll} plots the objective function with respect to the number of iterations independently for each algorithm and each database, using nested iteration as the multigrid cycle (which is the easiest to represent). 
In all cases, prolongations of solutions from the lower levels generate much better solutions than those obtained on the fine level (as explained in Section~\ref{smpro}).\\

These test results are very encouraging: the multilevel approach for NMF seems very efficient when dealing with image datasets and allows a significant speedup of the convergence of the algorithms.

\section{Limitations}

Although the numerical results reported in the previous section demonstrate significant computational advantages for our multilevel technique, we point out in this section two limitations that can potentially affect our approach. 

\subsection{Size of the Data Matrix} \label{sizeM}

The approach described above was applied to only one dimension of the input data: restriction and prolongation operators are applied to columns of the input matrix $M$ and of the first factor~$V$. Indeed, we assumed that each of these columns satisfies some kind of smoothness property. In contrast, we did not assume that the columns of $M$ are related to each other in any way, so that no such property holds for the rows of $M$. Therefore we did not apply our multilevel strategy along the second dimension of the input data, and our approach only reduced the row dimension $m$ of matrix $M$ at each level from $m$ to $m' \approx \frac{m}{4}$, while the column dimension $n$ remained the same. 

The fact that the row dimension of factor $V$ becomes smaller at deeper levels clearly implies that the computational cost associated with updating $V$ will decrease. This reduction is however not directly proportional to the reduction from $m$ to $m'$, as this cost also depends on the factorization rank $r$ and the dimensions of the other factor, which are not affected. Similarly, although the dimensions of factor $W$ remain the same regardless of the depth, its updates could become cheaper because dimension $m$ also plays a role there. The relative extent of those effects depends on the NMF algorithm used, and will determine in which situations a reduction in the dimension $m$ is clearly beneficial with respect to the whole computational cost of the algorithm. 

We now analyze in detail the effect of a reduction of $m$ on the computational cost of one iteration of the algorithms presented in Section~\ref{algonmf}: 
\begin{table}[ht!]
\caption{Number of floating point operations needed to update $V$ and $W$ in ANLS, MU and HALS.} 
\begin{center}
\renewcommand\arraystretch{1.2}
\begin{tabular}{c|c|c}
						  & ANLS       							&  MU and HALS   \\ \hline
Update of $V$	& $\mathcal{O}(mnr + m s(r) r^3 + n r^2)$     &   $\mathcal{O}(mnr + (m+n) r^2)$      			\\
Update of $W$	& $\mathcal{O}(mnr + n s(r) r^3 + m r^2)$     &   $\mathcal{O}(mnr + (m+n) r^2)$    \\
Both updates	& $\mathcal{O}( m (nr+s(r) r^3) + n s(r) r^3 )$     &   
$\mathcal{O}( m (nr + r^2) + n r^2 )$    \\
\end{tabular} 
 \\

\begin{flushleft}
\footnotesize (function $s(r)$ is $2^r$ in the worst case, and typically much smaller, see~\ref{AppA}).  
\end{flushleft}
\end{center}
\label{tableCC}
\end{table}  
\begin{figure}[ht!]
\begin{center}
\includegraphics*[width=.75\textwidth]{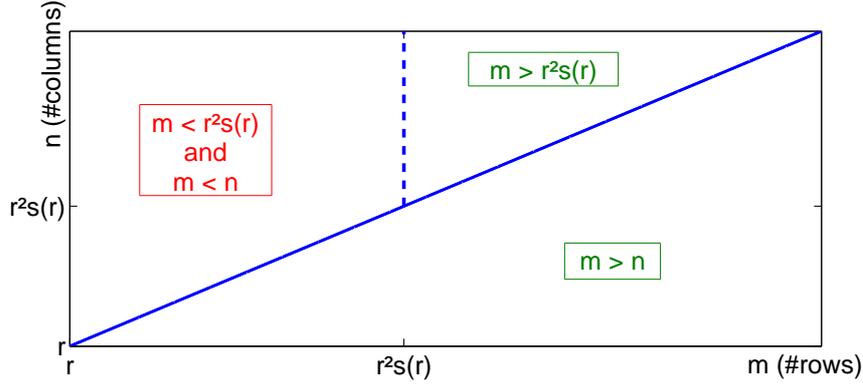}
\caption{Regions for input dimensions $(m,n)$ where a multilevel strategy is  beneficial in all cases ($m \geq \min\{n,r^2 s(r)\}$, lower and upper right parts) or only for MU and HALS ($m \leq \min\{n,r^2 s(r)\}$, upper left part).} 
\label{mnrZones} 
\end{center}
\end{figure}

Table~\ref{tableCC} gives the computational cost for the updates of $V$ and $W$ separately, as well as their combined cost (see~\ref{AppA}). Our objective is to determine for which dimensions $(m,n)$ of the input matrix and for which rank $r$ our multilevel strategy (applied only to the row dimension $m$) is clearly beneficial or, more precisely, find when a constant factor reduction in $m$, say $\frac{m}{m'}=4$, leads to a constant factor reduction in the total computational cost of both updates. 
We make the following observations, illustrated on Figure~\ref{mnrZones}. 

\begin{itemize}
\item We need only consider the region where both $m$ and $n$ are greater than the factorization rank $r$ (otherwise the trivial factorization with an identity matrix is optimal).

\item Looking at the last row of the table, we see that all terms appearing in the combined computational cost for both updates are proportional to $m$, except for two terms: $n s(r) r^3$ for ANLS and $n r^2$ for MU and HALS. If the contributions of those two terms could be neglected compared to the total cost, any constant factor reduction in dimension $m$ would lead to an equivalent reduction in the total complexity, which is the ideal situation for our multilevel strategy. 

\item When $m \ge n$, terms $n s(r) r^3$ for ANLS and $n r^2$ for MU and HALS are dominated respectively by $m s(r) r^3$ and $m r^2$ (i.e., $n s(r) r^3 \leq m s(r) r^3$ and $n r^2 \leq m r^2$), so that they cannot contribute more than half of the total computational cost. Therefore a reduction in dimension $m$ will guarantee a constant factor reduction in the total complexity. Let us illustrate this on the MU 
(a similar analysis holds for ANLS and HALS) 
for which the exact total computational cost is $2m(nr + r^2) + 2nr^2$ (see~\ref{AppA}). The factor reduction $f_{MU}$ in the total complexity satisfies  
\[
1 
\leq 
f_{MU} = \frac{m (nr + r^2) + n r^2}{m' (nr + r^2) + n r^2}
\leq \frac{m}{m'} = 4, 
\]  
and, for $m \ge n \ge r$ and $\frac{m}{m'}=4$, we have  that 
\[
f_{MU} 
\geq 
\frac{mnr + m r^2 + m r^2}{m'nr + m' r^2 + m r^2}
= 
\frac{4m' nr + 8m'  r^2 }{m'nr + 5m' r^2} 
\geq 
\frac{4m'  r^2 + 8m'  r^2 }{m' r^2 + 5m' r^2} 
=
 2, 
\] 
i.e., the total computational cost of the MU updates on the coarse level is at least twice cheaper than on the fine level.   
Moreover, when $m$ is much larger than $n$ ($m \gg n$), as is the case for our images, the terms in $n$ can be  neglected, and we find ourselves in the ideal situation described previously (with $f_{MU} \approx 4$). In conclusion, when $m \ge n$, we always have an appreciable reduction in the computational cost. 

\item Looking now at MU and HALS when $m$ is smaller than $n$, we see that the term $n r^2$ is always dominated by $m n r$ (i.e., $n r^2 \leq m n r$), because $m \ge r$ always holds. We conclude that a constant factor reduction in the total complexity can also be expected when $m$ is reduced. For example, for MU, we have 
\[
f_{MU}  
\geq 
\frac{mnr + m r^2 + mnr}{m'nr + m' r^2 + mnr} 
=
\frac{8m' nr + 4m' r^2 }{5m'nr + m' r^2} 
\geq \frac{8}{5}. 
\] 

\item Considering now ANLS when $m$ is smaller than $n$, we see that the term $n s(r) r^3$ is dominated by $m n r$ as soon as $m \ge s(r) r^2$. Again, in that situation, a constant factor reduction in the total complexity can be obtained\footnote{It is worth noting that when $m \geq s(r) r^2$ the initial computational cost to formulate the NNLS subproblem in $W$: 
\begin{equation} \label{subNNLS}
\min_{W \geq 0} \quad \sum_{i=1}^n ||M_{:i}-VW_{:i}||_F^2 = \sum_{i=1}^n ||M_{:i}||_F^2 - 2 (M_{:i}^T V)  W_{:i} + W_{:i}^T (V^T V) W_{:i}, 
\end{equation}
which requires the computation of $V^TV$ and $M^T V$ (cf.~\ref{AppA}), takes more time than actually solving \eqref{subNNLS}.}. Finally, the only situation where the improvement due to the multilevel technique is modest is when using ANLS when both $m < n$ and $m < s(r) r^2$ hold, in which case the term $n s(r) r^3$ can dominate all the others, and a reduction in dimension $m$ is not guaranteed to lead to an appreciable reduction in the total complexity. 

\end{itemize} 

To summarize, applying multilevel techniques to the methods presented in this paper is particularly
beneficial on datasets for which $m$ is sufficiently large compared to $n$ and $r$ (for MU and HALS) and to $n$ and $s(r) r^2$ (for ANLS). Some gains can always be expected for MU and HALS, while ANLS will only see a significant improvement if $m \ge \min\{n, s(r) r^2\}$ holds.

In Section~\ref{Res}, we have  presented computational experiments for image datasets satisfying this requirement: the number of images $n$ was much smaller than the number of pixels $m$ in each image. In particular, we observed that the acceleration provided by the multilevel approach to the ANLS algorithm was not as significant as for HALS: while in most cases ANLS converged faster than HALS when using the original NMF algorithms, it converged slower as soon as the multilevel strategy was used (see Tables~\ref{errortime1} and \ref{errortime3}).\\

To end this section, we note that, in some applications, rows of matrix $M$ can also be restricted to lower-dimensional spaces. In these cases, the multilevel method could be made even more effective. This is the case for example in the following situations:  
\begin{itemize}

\item In hyperspectral images, each column of matrix $M$ represents an image at a given wavelength, while each row represents the spectral signature of a pixel, see, e.g., \cite{hubble1, GP11}. Since spectral signatures feature smooth components, the multilevel strategy can be easily generalized to reduce the number of rows $n$ of the data matrix $M$. 

\item For a video sequence, each column of matrix $M$ represents an image at a given time so that consecutive images share similarities. Moreover, if the camera is fixed, the background of the scene is the same among all images. The multilevel approach can then also be generalized to reduce the number of columns of $M$ in a meaningful way.  

\item In face datasets (e.g., used for face recognition), a person typically appears several times.  Hence one can imagine using the multilevel strategy by merging different columns corresponding to the same person. 

\end{itemize}

%

\subsection{Convergence} \label{conv}

In classical multigrid methods, when solving a linear system of equations $Ax = b$,  the current approximate solution $x_c$ is not transferred from a fine level to a coarser one, because it would imply the loss of its high-frequency components; instead, the residual is transferred, which we briefly explain here. 
Defining the current residual $r_c = b-Ax_c$ and the error $e = x-x_c$, we have the equivalent defect equation $Ae = r_c$ and we would like to approximate $e$ with a correction $e_c$ in order to improve the current solution with $x_c \leftarrow x_c+e_c$. The defect equation is solved approximately on the the coarser grid by restricting the residual $r_c$, 
the correction obtained on the coarser grid is prolongated and the new approximation $x_c+e_c$ is computed, see, e.g., \cite[p.37]{Trot}. 
If instead the solution is transferred directly from one level to another (as we do in this paper), the corresponding scheme is in general not convergent, see \cite[p.156]{Trot}. In fact, even an exact solution of the system $Ax=b$  is not a fixed point, because the restriction of $x$ is not an exact solution anymore at the coarser level (while, in that case, the residual $r$ is equal to zero and the correction $e$ will also be equal to zero). 

Therefore, the method presented in this paper should in principle only be used as a \emph{pre-processing} or \emph{initialization} step before another (convergent) NMF algorithm is applied. 
In fact, if one already has a good approximate solution $(V,W)$ for NMF (e.g.,  a solution close to a stationary point),
 transferring it to a coarser grid will most likely increase the approximation error  because  high frequency components (such as edges in images) will be lost.  Moreover, it seems that the strategy  of transferring a residual instead of the whole solution is not directly applicable to NMF. Indeed, a  `local linearization' approach, which would consist in  linearizing the equation
\[ 
M-(V+\Delta V)(W+\Delta W) \approx 0 \quad \iff \quad R = M-VW \approx V \Delta W+\Delta V W, 
\] 
where $\Delta V$ and $\Delta W$ are the corrections to be computed on the coarser grids, causes several problems. First, handling non-negativity of the coarse versions of the factors becomes non-trivial. Second, performing this approximation efficiently also becomes an issue, since for example computing the residual $R$ is as expensive as computing directly a full MU or HALS iteration on the fine grid ($\mathcal{O}(mnr)$ operations). Attempting to fix these drawbacks, which seems to be far from trivial, is a topic for further research.  

To conclude this section, we reiterate that, despite these theoretical reservations, it seems our technique is still quite efficient (see Section~\ref{Res}). One reason that explains that good behavior is that NMF solutions are typically \emph{part-based} and \emph{sparse} \cite{LS1}, see Figure~\ref{exNMF}.  Therefore, columns of matrix $V$ contains relatively large `constant components', made of their zero entries, which are perfectly transferred from one level to another, so that $s_V = \frac{|| V - \mathcal{P}(\mathcal{R}(V))||_F}{||V||_F}$ will typically be very small 
(in general much smaller than $s_M$). 

\section{Concluding Remarks} 

In this paper, a multilevel approach designed to accelerate NMF algorithms has been proposed and its efficiency has been experimentally demonstrated. Applicability of this technique relies on the ability to design linear operators preserving nonnegativity and transferring accurately data between different levels. To conclude, we give some directions for further research.

\subsection{Extensions} \label{ext}

We have only used our multilevel approach for a specific objective function (sum of squared errors) to speed up three NMF algorithms (ANLS, MU and HALS) and to factorize 2D images. However, this technique can be easily generalized to different objective functions, other iterative algorithms and applied to various kinds of smooth data. In fact, the key characteristic we exploit is the fact that a reduction of the dimension(s) of the input matrix (in our numerical examples, $m$) leads to cheaper iterations (on coarse levels) for any reasonable algorithm, i.e., any algorithm whose computational cost depends on the dimension of the input matrix (see also the more detailed analysis in Section~\ref{sizeM}). 

Moreover, other types of coarse grid definition (e.g., red-black distribution), 
transfer operators (e.g., wavelets transform) and grid cycle (e.g., W-cycle or flexible cycle) 
can be used and could potentially further improve efficiency. 

This idea can also be extended to nonnegative tensor factorization (NTF), see, \rm{e.g.}, \cite{CAZP09, hubble} and the references therein, by using multilevel techniques for higher dimensional spaces. 

\subsection{Initialization}

Several judicious initializations for NMF algorithms have been proposed in the literature which allow to accelerate convergence and, in general, improve the final solution~\cite{Cur,  Bou}. 
The computational cost of these good initial guesses depends on the matrix dimensions and will then be cheaper on a coarser grid. Therefore, it would be interesting to combine classical NMF initializations techniques with our multilevel approach for further speedups. 

\subsection{Unstructured data}

When we do not possess any kind information about the matrix to factorize (and a fortiori about the solution), applying a multilevel method seems out of reach. In fact, in these circumstances, there is no sensible way to define the transfer operators.

Nevertheless, we believe it is not hopeless to extend the multilevel idea to other types of data. For example, in text mining applications, the term-by-document matrix can be restricted by stacking synonyms or similar texts together, see \cite{SFS09} where graph coarsening is used. This implies some a priori knowledge or preprocessing of the data and, assuming it is cheap enough, the application of a multilevel strategy could be expected to be profitable in that setting. 

\subsection*{Acknowledgments}

We thank Quentin Rentmeesters, Stephen Vavasis and an anonymous reviewer for their insightful comments which helped improve the paper.

\bibliographystyle{siam}
\bibliography{multilevel}

\normalsize


\appendix

\section{Computational Cost of ANLS, MU and HALS} \label{AppA}

\subsection{MU and HALS}

The main computational cost for updating $V$ in both MU and HALS resides in the computation of $MW^T$ and\footnote{In HALS, $VWW^T$ is essentially computed one column at a time, see \cite[Section 4.2.1]{Gillis}.} $V(WW^T)$,  which requires respectively $2mnr$ and $2(m+n)r^2$ operations, cf.\@ Algorithms~\ref{MU} and \ref{HALS}. Updating $W$ requires the same number of operations, so that the total computational cost is $\mathcal{O}(mnr+(m+n)r^2)$ operations per iteration, almost proportional to $m$ (only the $nr^2$ term is not, but is negligible compared to the other terms, cf.\@ Subsection~\ref{sizeM}), see also \cite[Section 4.2.1]{Gillis}. 

\subsection{Active-Set Methods for NNLS}

In a nutshell, active-set methods for nonnegative least squares work in the following iterative fashion~\cite[Algorithm NNLS, p.\@ 161]{LH74}
\begin{enumerate}
\item[0.] Choose the set of active (zero) and passive (nonzero) variables. 

\item Get rid of the nonnegativity constraints and solve the unconstrained least squares problem (LS) corresponding to the set of passive (nonzero) variables (the solution is obtained by solving a linear system, i.e., the normal equations); 

\item Check the optimality conditions, i.e., the nonnegativity of passive variables, and the nonnegativity of the gradients of the active variables. If they are not satisfied: 

\item Exchange variables between the set of active  and the set of passive variables in such a way that the objective function is decreased at each step; and go to 1. 
\end{enumerate} 

In \eqref{NMF}, the problem of computing the optimal $V$ for a given fixed $W$ can be decoupled into $m$ independent NNLS subproblems in $r$ variables:
\[
\min_{V_{i:} \in \mathbb{R}^r_+} ||M_{i:}-V_{i:}W||_F^2, \; 1 \leq i \leq m. 
\]
Each of them amounts to solving a sequence of linear subsystems (with at most $r$ variables, cf.\@ step 1 above) of 
\[
V_{i:} (WW^T) = M_{i:}W^T, \; 1 \leq i \leq m. 
\] 
In the worst case, one might have to solve every possible subsystem, which requires $\mathcal{O}(g(r))$ operations with\footnote{One can check that $(2^{(r-3)}-1) (r-2)^3 \leq g(r) \leq 2^r r^3$.}
$g(r) = \sum_{i=1}^r \binom{r}{i} i^3$ = $\Theta(2^r r^3)$. 
Note that $WW^T$ and $MW^T$ 
can be computed once for all, which requires $\mathcal{O}(mnr + nr^2)$ operations (see previous section on MU and HALS). Updating $V$ then requires $\mathcal{O}(mnr + m s(r)  r^3 + nr^2)$ operations, while updating $W$ similarly requires $\mathcal{O}(mnr + n s(r)  r^3 +  mr^2)$. Finally, the total computational cost of one ANLS step is $\mathcal{O}(mnr + (m+n) r^2 (r s(r) +1)) = \mathcal{O}(mnr + (m+n) s(r)  r^3)$ operations per iteration, where $s(r) \leq 2^r$. The number of steps $s(r)$ is $\Theta(2^r)$ in the worst case, but in practice is typically much smaller (as is the case for the simplex method for linear programming). \\

When $m$ is reduced by a certain factor (e.g., four as in the multilevel approach presented in Section~\ref{multiNMF}), the computational cost is not exactly reduced by the same factor, because the leading $(m+n)$ factor above also depends on $n$. However, in our applications, when $m$ (number of pixels) is much larger than $n$ (number of images), one can roughly  consider the cost per iteration to be reduced by the same factor, since $\frac{m+n}{4} \approx \frac{m}{4}$ (see also the more detailed discussion in Subsection~\ref{sizeM}).

\end{document}